\documentstyle{article}

\newtheorem{defn1}{ \sc Definition}[section]
\newtheorem{exmp1}[defn1]{ \it Example}
\newtheorem{prop1}[defn1]{ \sc Proposition}
\newtheorem{thm1}[defn1]{ \sc Theorem}
\newtheorem{cor1}[defn1]{ \sc Corollary}
\newtheorem{rem1}[defn1]{ \it Remark}

\def\mpr#1{\;\smash{\mathop{\hbox to 20pt{\rightarrowfill}}\limits^{#1}}\;}
\def\mpd#1{\big\downarrow\rlap{$\vcenter{\hbox{$\scriptstyle#1$}}$}}
\def\mdp#1{\llap{$\vcenter{\hbox{$\scriptstyle#1$}}$}\big\downarrow}
\def\mup#1{\llap{$\vcenter{\hbox{$\scriptstyle#1$}}$}\big\uparrow}

\def\D{\,I\!\!\!\!D}
\def\IM{I\!M}
\def\IC{I\!C}
\def\IH{I\!H}
\def\K{K\!E\!R}
\def\o{{\rm reg}}

\def\bar{\overline}
\def\p{\phantom{\mpd{}^\mpd{}_\mpd{}}}
\def\qed{\hfill$\Box$\vskip10pt}

\def\Cset{{\bf C}}
\def\Pset{{\bf P}}
\def\Zset{{\bf Z}}
\def\Rset{{\bf R}}
\def\Nset{{\bf N}}
\def\Qset{{\bf Q}}

\begin{document}

\author{Andrzej Weber\footnote{Supported by KBN 2P03A 00218 grant.
I thank Institute of Mathematics, Polish Academy of Science for
hospitality.}\\
\small\it Instytut Matematyki, Uniwersytet Warszawski }
\date{}
\title{Pure homology of algebraic varieties}
\maketitle

\begin{abstract} We show that for a complete complex algebraic
variety the pure term of the weight filtration in homology 
coincides with the image of intersection homology. Therefore pure
homology is topologically invariant. To obtain slightly more general
results we introduce {\it image homology} for noncomplete varieties.

\vskip 6pt
\noindent {\sc Key words}: Algebraic varieties, weight filtration,
intersection homology.
\end{abstract}

\section{Introduction}

The mixed Hodge theory was developed by P.~Deligne to study
extraordinary properties of cohomology of complex algebraic varieties. One
of the ingredients of mixed Hodge structure is the weight
filtration. We will focus on the dual filtration in homology.  We
will find a relation between intersection homology defined by
M.~Goresky--R.~MacPherson and the weight filtration. The main
result of the paper says that if a variety is complete, then the
pure part of homology $W^kH_k(X)$ is the image of intersection
homology. This shows topological invariance of pure homology.  An
attempt to give a topological description or at least estimate the
weight filtration was undertaken by C.~McCrory (\cite{MC}). If $X$
is a complete smooth divisor with normal crossings then the weight
filtration coincides with the Zeeman filtration given by
"codimension of cycles". In general there is an inclusion (the
terms of the weight filtration are bigger). It is a result of McCrory
for hypersurfaces and of F.~Guillen (\cite{Gu}) in general. For
complete, normal and equidimensional varieties all we can deduce
from Zeeman filtration for the pure term of the weight filtration
is the following: the image of the Poincar\'e duality map is
contained in $W^kH_k(X)$. This is exactly what follows from the
weight principle. Our description of $W^kH_k(X)$ is complete and
still purely topological. Other terms of the weight filtration do
not have such description. They are not topologically invariant.

In the course of the proof at some points we do not have to assume
that the variety is complete. Therefore we organize the paper as
follows.  Let $X$ be a complex algebraic variety. We will define
certain subspaces of homology
$$\IM_k(X)\subset H_k(X)$$ for each $k\in\Nset$, which we will call
{\it image homology}. We assume that
homology has coefficients
in a field of characteristic zero. We will show that:
\begin{thm1}\label{th1} The image homology satisfies the 
conditions:\begin{description}
\item[\sc (a)] If $X$ is smooth, then $\IM_k(X)=H_k(X)$, 
\item[\sc (b)] If $f:X\rightarrow Y$ is an algebraic map then
$f_*\IM_k(X)\subset \IM_k(Y)$,
\item[\sc (c)] If $f:X\rightarrow Y$ is an algebraic map which is proper and
surjective then $f_*\IM_k(X)=\IM_k(Y)$.
\end{description}\end{thm1}
By Hironaka resolution of singularities every algebraic variety is dominated
by a smooth one.
Therefore the conditions of Theorem \ref{th1} determine
$\IM_k(X)$. Of course it is not clear whether $\IM_k(X)$
satisfying {\sc(a)-(c)} exists. If $X$ and $Y$ are complete then
the existence of such subspaces is a consequences of mixed Hodge
theory, \cite{De} \S8.2.
\begin{thm1}\label{th2}
\begin{description}\item[]
\item[\sc (d)] If $X$ is complete, then $\IM_k(X)$ coincides with
$W^kH_k(X)$, the pure weight subspace of homology.
\end{description}\end{thm1}
We prove that $\IM_k(X)$ is the image of intersection homology:
\begin{thm1}\label{th3}
\begin{description}\item[]
\item[\sc (e)] If $X$ is equidimensional, then
$$\IM_k(X)=im(\iota_X:\IH_k(X)\rightarrow H_k(X))\,.$$
\end{description}\end{thm1}
We will explain the objects appearing in {\sc (d)} and {\sc (e)} below:
\begin{description}
\item [\rm (d)] $W^kH_k(X)=W_{-k}H_k(X)$ is the $k$-th term of the
weight filtration in homology of $X$. It is the annihilator of
$W_{k-1}H^k(X)$,
the subspace of cohomology of the weight $\leq k-1$,
defined in \cite{De}. If $X$ is
complete then $W^kH_k(X)$ is the lowest possibly nonzero term. Then
it is pure. In
general the weight filtration is not a topologically invariant, see \cite{St}.
\item [\rm (e)] $\IH_k(X)$ is the intersection (co)homology group
(with respect to the middle perversity)
defined in \cite{GM1,GM2} and developed in \cite{BBD}. It is a
cohomology theory adapted to tackle singular varieties. It is
equipped with canonical maps:
$$H^{2\dim(X)-k}_c(X)\mpr{\iota_X^*}\IH^{2\dim(X)-k}_c(X)
=\IH_k(X)\mpr{\iota_X}H_k(X)$$
factorizing Poincar\'e duality map $[X]\cap-$.
(The subscript $c$ stands for compact supports.)
Intersection (co)homology is a topological invariant of
$X$.
According to \cite{BBD} or \cite{Sa} one can construct a
weight filtration in $\IH_k(X)$. If $X$ is complete, then $\IH_k(X)$ is
pure of weight $k$.
\end{description}
We deduce a striking property of
image homology.
\begin{thm1}[Topological Invariance] \label{tht}
\begin{description}\item[ ]
\item[\sc (f)] The image homology is a topological
invariant of $X$.
This means, that if $f:X\rightarrow Y$ is a homeomorphism which does
not have to be an algebraic map then $f_*\IM_k(X)=\IM_k(Y)$.
\end{description}\end{thm1}
Theorem \ref{tht} follows from {\sc (e)}. As a corollary we obtain an unusual
property of the pure term of
the weight filtration:
\begin{cor1} If X is a complete algebraic variety, then $W^kH_k(X)$
is topologically invariant.\end{cor1} A topological description of
$W^kH_k(X)$ is the following
\begin{cor1} If X is a complete and equidimensional algebraic variety,
then
$$W^kH_k(X)=im(\iota_X:\IH_k(X)\rightarrow H_k(X))\,.$$
\end{cor1}

We also define parallel or dual variants of $\IM_k(X)$:
\begin{description}
\item [\rm --] $\IM^{BM}_k(X)\subset H^{BM}_k(X)$, a subspace of Borel-Moore
homology satisfying the properties {\sc (a)-(f)} (with {\sc (b)}
for proper maps),
\item [\rm --] $\K^k(X)\subset H^k(X)$, {\it kernel cohomology}, a subspace
of cohomology satisfying dual properties,
\item [\rm --] $\K^k_c(X)\subset H_c^k(X)$, a subspace of cohomology with
compact supports satisfying dual properties
(with {\sc (b')} for proper maps).
\end{description}
The dual properties are the following:
\begin{thm1}\label{th4} The kernel cohomology $\K^k(X)\subset
H^k(X)$ satisfies the conditions:
\begin{description}
\item[\sc (a')] If $X$ is smooth, then $\K^k(X)=0$,
\item[\sc (b')] If $f:X\rightarrow Y$ is an algebraic map then
$f^*\K^k(Y)\subset \K^k(X)$,
\item[\sc (c')] If $f:X\rightarrow Y$ is an algebraic map which is proper and
surjective then $(f^*)^{-1}\K^k(X)=\K^k(Y)$,
\item[\sc (d')] If $X$ is complete, then $\K^k(X)=W_{k-1}H^k(X)$,
the subspace of cohomology of the weight $\leq k-1$,
\item[\sc (e')] If $X$ is equidimensional, then
$$\K^k(X)=ker(\iota_X^*:H^k(X)\rightarrow \IH^k(X))\,,$$
\item[\sc (f')] The kernel cohomology is a topological invariant of $X$.
\end{description}\end{thm1}
Instead of $\K$ groups one can consider the image of cohomology
in intersection
cohomology. We will focus on the case when $X$ is equidimensional.
The property {\sc (d')} can be restated:
\begin{thm1}
\begin{description}
\item[\ ]
\item[\sc (d'')] If $X$ is complete and equidimensional then the
image $im(\iota^*_X:H^k(X)\rightarrow\IH^k(X))$ is isomorphic to the
pure quotient $H^k(X)/W_{k-1}H^k(X)$.
\end{description}\end{thm1}
Combining {\sc(a')} and {\sc(c')} we obtain
\begin{cor1}\label{co} If $X$ is smooth, $Y$ is equidimensional and the map
$f:X\rightarrow Y$ is proper and surjective, then the image
$f^*(H^*(Y))$ is isomorphic to the image $\iota_X^*$.\end{cor1}
Corollary \ref{co} (and  an analogous one for homology) allows us
to think about intersection cohomology as a substitute for the
cohomology of a minimal resolution.

The main tool of the proofs is the Decomposition Theorem of \cite{BBD}
or \cite{Sa}. We will rely only on one corollary:

\begin{cor1}[from the Decomposition Theorem]\label{cdt} Let
$f:X \rightarrow Y$ be a proper surjective map of algebraic
varieties. Then $\IC_Y$ is a retract of $Rf_*\IC_X$, i.e.~there
exist maps $i:\IC_Y\rightarrow Rf_*\IC_X$ (inclusion) and
$r:Rf_*\IC_X\rightarrow \IC_Y$ (retraction) such that $r\circ
i=Id_{\IC_Y}$.
\end{cor1}
Here $\IC_X$ is the intersection complex. It is an object of the derived
category $D(X)$ of sheaves of vector spaces over $X$ and $X$ is assumed to
have pure dimension $n$. The hypercohomology of $\IC_X$ is the
intersection (co)homology:
$$H^k(X;\IC_X)=\IH^k(X)=\IH^{BM}_{2n-k}(X)\,,$$
$$H_c^k(X;\IC_X)=\IH_c^k(X)=\IH_{2n-k}(X)\,.$$
Due to methods of \cite{BBD} we are forced to use complex
coefficients (or $\bar\Qset_\ell$), but once we prove the Theorems
\ref{th1}-\ref{tht} for $\Cset$ they will automatically follow for
arbitrary field of characteristic zero. We will also apply a
result of \cite{BBFGK} concerning (non canonical) functoriality of
intersection (co)homology. As it is shown in \cite{We1} the result
of \cite{BBFGK} is a formal consequence of \ref{cdt}. On the other
hand Corollary \ref{cdt} can be deduced from \cite{BBFGK}. We
present this reasoning in the Appendix. This is the only
non topological ingredient we use. The methods presented here are
applied to study residues on singular hypersurfaces in \cite{We3}.

\section{Definition through {\sc(e)}}
The starting point of our definition is the property {\sc (e)}. We
cannot say that $\IM_k(X)$ is just the image of $\IH_k(X)$ since
intersection homology is defined only for equidimensional
varieties. Therefore we decompose $X$ into irreducible components
$X=\bigcup_{i\in I} X_i$.
\begin{defn1} Let $\iota_i:\IH_k(X_i)\rightarrow H_k(X)$
be the composition of the natural transformation $\iota_{X_i}$
with the map induced by the inclusion of the component. The image
homology of $X$ is defined by
$$\IM_k(X)=\sum_{i\in I}im(\iota_i)\subset H_k(X)\,.$$
\end{defn1}
If $X$ is equidimensional we do not have to decompose $X$ into
components since
$$\IH_k(X)=\IH_k(\widehat X)=\bigoplus_{i\in I}\IH_k(X_i)\,,$$
where $\widehat X$ is the normalization of $X$, which is also the
normalization of $\bigsqcup_{i\in I}X_i$. This proves the property
{\sc (e)}.

We can give a direct description of image homology in terms of
cycles. Our description is based on the elementary definition of
intersection homology introduced in \cite{GM1}. We fix a
stratification $\{S_\alpha\}_{\alpha\in J_i}$ of each $X_i$. We
assume that the stratification is locally topologically trivial. A
homology class of $H_k(X)$ belongs to $\IM_k(X)$ if it is
represented by a cycle $\xi$ which
\begin{description}
\item[--] can be decomposed into a sum of cycles $\xi=\sum_{i\in
I}\xi_i$, with $|\xi_i|\subset X_i$.
\item[--] the component $\xi_i$ is allowable in $X_i$, i.e.
$$\dim_\Rset(|\xi_i|\cap S_\alpha)\leq k-c-1$$
if $S_\alpha$ is a singular stratum of $X_i$ with complex
codimension $c=\dim(S_\alpha)-\dim(X_i)$ .
\end{description}

\section{Proof of \sc(f)}
By \cite{GM2} or \cite{Ki} intersection homology is a topological
invariant. It is enough to show that the decomposition into
irreducible components is topologically invariant. Indeed, let
$X_\o$ be the set of points of $X$ at which $X$ is a topological
manifold (the dimension may vary from a point to a point). The set
$X_\o$ decomposes into a set of components $X_{i,\o}$. The
irreducible components of $X$ are the closures of $X_{i,\o}$.\qed

\section{Proof of \sc(a)}
If $X$ is smooth, then $\iota_X:\IH_k(X)\rightarrow H_k(X)$ is an
isomorphism. \qed
\begin{rem1}\rm It is enough to assume that $X$ is a rational homology
manifold to have $\IM_k(X)=H_k(X)$. Also some information about
local homology allows to deduce an isomorphism in certain degrees.
For example if $X$ has isolated singularities, then
$\IM_k(X)=H_k(X)$ for $k>\max_i(\dim(X_i))$.\end{rem1}

\section{Proof of \sc(b)}
At this point for the first time we have to use highly nontrivial
results concerning intersection homology. Let $f:X\rightarrow Y$
be an algebraic map. Let us recall the main theorem of \cite{BBFGK}.
\begin{thm1}\label{func}There is a map $\beta$ making the following
diagram commute: 
$$\matrix{\IH_k(X)&\mpr{\iota_X}&H_k(X)\cr
\mdp{\beta}&\p&\mpd{f_*}\cr
\IH_k(Y)&\mpr{\iota_Y}&H_k(Y)&.\cr}$$\end{thm1} 
The choice of $\beta$ is not canonical, but its existence is enough
to deduce {\sc(b)}.\qed

In some cases we can prove more:

\begin{prop1} \label{spw} If $X$ is complete, then $f_*$ strictly
preserves image homology:
$$im(f_*)\cap \IM_k(Y)=f_*(\IM_k(X))\,.$$
\end{prop1}
We will deduce \ref{spw} from mixed Hodge theory of \cite{De}
provided we know {\sc(c)} and {\sc(d)}.

{\noindent\hskip 6pt \it Proof.} Let $\mu:\widetilde Y\rightarrow Y$ be a
resolution of singularities. By {\sc(c)} $\IM_k(Y)$ is the image
of $\mu_*$. Therefore  $\IM_k(Y)\subset W^kH_k(Y)$. The map $f_*$
strictly preserves weights
$$im(f_*)\cap W^kH_k(Y)=f_*(W^kH_k(X))\,.$$
By {\sc(d)} we have $W^kH_k(X)=\IM_k(X)$. We obtain the inclusion
$$im(f_*)\cap \IM_k(Y)\subset f_*(\IM_k(X))\,.$$
The converse inclusion follows from {\sc(b)}. \qed

\section{Proof of \sc(c)}\label{szesc}
We assume that $f:X\rightarrow Y$ is proper and surjective.
It is enough to prove {\sc(c)} for irreducible $X$ and $Y$.
We will show that
the map $\iota_Y:\IH_k(Y)\rightarrow H_k(Y)$ factors through $IH_k(X)$.
Our argument is dual to the one of \cite{We1}.
We will construct a map $\alpha$ which will fit to the commutative
diagram
$$\matrix{\IH_k(X)&\mpr{\iota_X}&H_k(X)\cr
\mup{\alpha}&\p&\mpd{f_*}\cr
\IH_k(Y)&\mpr{\iota_Y}&H_k(Y)&.\cr}$$
It will come from a map $\bar\alpha$ at the level of derived category
of sheaves:
$$\matrix{Rf_*\IC_X[2n]&\mpr{Rf_*(\iota_X)}&Rf_*\D_X\cr
\mup{\bar\alpha}&\p&\mpd{f_*}\cr
\IC_Y[2m]&\mpr{\iota_Y}&\D_Y&.\cr}$$ Here $n=\dim(X)$, $m=\dim(Y)$
and $\D_X$ (resp. $\D_Y$) denotes the dualizing sheaf. It is equal
to $\Cset_X[2n]$ (resp. $\Cset_Y[2m]$) at smooth points. By the
Corollary from the Decomposition Theorem \ref{cdt} there is a
retraction $r:Rf_*\IC_X\rightarrow\IC_Y$. Let
$$\bar\alpha:\IC_Y[2m]=D\IC_Y\mpr{Dr}
DRf_*\IC_X=Rf_*D\IC_X=Rf_*\IC_X[2n]$$ be the Verdier dual of the
retraction. We obtain two maps: $\iota_Y$ and the composition $f_*\circ
Rf_*(\iota_X)\circ \bar\alpha:\IC_Y[2m]\rightarrow\D_Y$. We can
rescale  $\bar\alpha$ and assume that these two maps coincide on
an open set. We will show that they coincide on the whole $X$. The
conclusion will follow.

\begin{prop1}\label{ro} Suppose $Y$ is irreducible. Then every two maps
in the derived category $\iota_1,\iota_2:\IC_Y\rightarrow\D_Y$
which coincide on an open set are equal.\end{prop1}

{\noindent\hskip 6pt \it Proof.} It is more convenient to prove the dual statement: {\it 
Every
two maps $D\iota_1,D\iota_2:\Cset_Y\rightarrow\IC_Y$ which
coincide on an open set are equal.} Indeed
$${\rm Hom}_{D(Y)}(\Cset_Y,\IC_Y)=R^0Hom(\Cset_Y,\IC_Y)
=H^0(Y;\IC_Y)=\IH^0(Y)\simeq\Cset$$
and the same for $U$. The restriction map $\IH^0(Y)\rightarrow
IH^0(U)$ is an isomorphism.\qed

\begin{rem1}\rm Another possible proof is to show that $\beta$ in
\ref{func} may be chosen to be surjective.\end{rem1}

\begin{rem1}\label{sup}\rm Since the factorization of $\iota_Y$ was 
constructed
on
the sheaf level we obtain the statement for arbitrary supports of the
homology:
$$f_*\IM_k^{f^{-1}(\phi)}(X)=\IM_k^\phi(Y)\subset H_k^{\phi}(Y)\,,$$
where $\phi$ is a family of supports in $Y$.
In particular
$$f_*\IM_k^{BM}(X)=\IM_k^{BM}(Y)\subset H_k^{BM}(Y)\,,$$
provided that the map $f$ is proper.
\end{rem1}

\begin{rem1}\rm\label{ac} Every algebraic subvariety $A\subset Y$ of dimension
$k$ defines a class in $\IM^{BM}_{2k}(Y)$. It is easy to show a
class in $\IM^{BM}_{2k}(X)$ which is mapped to $[A]$. To construct
it we choose an open and dense subset $U\subset A$ and a
subvariety $B\subset f^{-1}U$ such that $f_{|B}:B\rightarrow U$ is
finite, say of degree $d$. Let $C$ be the closure of $B$ in $X$.
Then the class $\frac1d [C]\in \IM^{BM}_{2k}(X)$ is the desired
lift of $[A]$.\end{rem1}

\begin{rem1}\rm\label{node} The property {\sc (c)} does not hold for the maps
which are not algebraic. Example: Let $N$ be the rational node,
i.e. $\Pset^1$ with two points identified. It is homeomorphic with
the famous pinched torus, i.e. $S^1\times S^1$ with $S^1\times pt$
shrunk. The quotient map $S^1\times S^1\rightarrow N$ is
surjective on $H_1(N)\simeq\Cset$, but $\IM_1(N)=H_1(\Pset^1)=0$.
\end{rem1}

Let us note that the property {\sc (c)} does not hold in the
analytic category, even in the smooth case.

\begin{exmp1}\label{anal} \rm Let $X=(\Cset^2\setminus \{0\})/\Zset$, where 
$k\in \Zset$ acts
on $\Cset^2$ via the multiplication by $2^k$. The tautological
bundle factors through $f:X\rightarrow \Pset^2=Y$. The map $f$ is
a complex analytic locally trivial fibration with fibers
$\Cset^*/\Zset$ which is an elliptic curve (a compact torus).
Topologically $f:X\simeq S^3\times S^1 \rightarrow S^2\simeq Y$ is
the projection on the first factor composed with the Hopf
fibration. The induced map is not surjective on $H_2(Y)$.
\end{exmp1}

\section{Proof of \sc(d)}
We recall that the weight filtration in the homology of a complete
variety is such, that
$W^{k+1}H_k(X)=0$ and $W^kH_k(X)$ is pure of weight $k$.
Our proof is based on the following description of $W^kH_k(X)$ (see
\cite{MC} p.218): Assume that $X$ is contained in a smooth variety
$M$. Let
$\mu:\widetilde M \rightarrow M$ be a resolution of $(M,X)$ in the sense
that $\mu^{-1}X$ is a smooth divisor with normal crossings. Let $\bar
X$ be the disjoint union of the components of $\mu^{-1}X$. Then
we have an equality
$$W^kH_k(X)=\mu_*(W^kH_k(\mu^{-1}X))=\nu_*(H_k(\bar X))\,,$$
where $\nu:\bar X\rightarrow X$ is the obvious map.
On the other hand by {\sc(a)} and {\sc(c)}
$$\nu_*(H_k(\bar X))=\IM_k(X)\,.$$
\qed

\begin{rem1}\rm By the dual of \cite{De}, 8.2.5 we can replace
$\nu:\bar X\rightarrow X$
by any surjective map from a complete and smooth variety.
\end{rem1}

\begin{rem1}\rm In general, if $X$ is possibly not complete then
$$\IM_k(X)\subset W^kH_k(X)\,.$$ 
The equality does not have to hold.
Example: Let $X=\Cset^*\times N$, where $N$ is the rational node
(see \ref{node}). Then $W^2H_2(X)=H_2(X)\simeq\Cset^2$, but
$\IM_2(X)\simeq H_2(\Cset^*\times\Pset^1)\simeq\Cset$.\end{rem1}

\section{Borel-Moore homology}
Proof of the properties {\sc(a)}, {\sc(d)}, {\sc(e)} and {\sc (f)}
are the same.
Concerning {\sc(b)} and {\sc(c)}, as we have remarked in \ref{sup}, we
can use any kind of supports, since the maps are constructed on the
level of sheaves.

\section{Dual $\K$ groups}
Intersection cohomology is a module over $H^*(X)$.
Since the natural transformation $\iota_X$ preserves the
module structure we have:
\begin{prop1} The kernel $ker(\iota_X)=\K^*(X)$ is an ideal in the
cohomology ring $H^*(X)$. Moreover $\K^0(X)=0$.\end{prop1} Due to
nondegeneracy of the pairing in intersection (co)homology we can
describe kernel cohomology as annihilators
$$\K^k(X)=Ann(\IM_k(X))\,,\qquad \K^k_c(X)=Ann(\IM_k^{BM}(X))\,.$$
Properties of {\sc (a')-(f')} can be obtained by duality.
A short way to prove {\sc(c')} is to show that $\beta$ in the
dual of the diagram \ref{func} may be chosen to be injective.
The inclusion of \ref{cdt} is a good choice.
We will comment on {\sc (d')}:
We recall that the weight filtration in the cohomology of a complete
variety is such, that
$W_kH^k(X)=H^k(X)$ and the quotient $W_kH^k(X)/W_{k-1}H^k(X)$ is pure.
Combining {\sc(d')} and {\sc(e')} we have:

\begin{thm1} If $X$ is complete and equidimensional, then
$$W_{k-1}H^k(X)=ker(\iota_X^*:H^k(X)\rightarrow IH^k(X))$$\end{thm1}
This gives a positive answer to a question stated in \cite{CGM}.

\section{Differential}
By the property {\sc(b)} the image homology is functorial. (We
recall that it is not possible to chose maps of intersection
homology to make it functor, see \cite{We1}.) Moreover the map
$\oplus\iota_i:\bigoplus_{i\in I}\IH_k(X_i)\rightarrow H_k(X)$
comes from a map of sheaves $\oplus\iota_i:\bigoplus_{i\in
I}\IC_{X_i}\rightarrow \D_X$. Therefore
$\partial(\IM_k(X))\subset\IM_{k-1}(A\cap B)$ in the
Mayer-Vietoris sequence associated to the open covering $X=A\cup
B$. Of course the sequence of $\IM$ groups does not have to be
exact. Example: Let $X$ the sum of two copies of $\Pset^1$ glued
along two pairs points. Denote by $a$ and $b$ be the points of the
identification. The Mayer-Vietoris sequence for
$A=X\setminus\{b\}$, $B=X\setminus\{a\}$ is not exact at
$\IM_0(A\cap B)$.

On the other hand one can consider the case when $A$ and $B$ are
closed subvarieties. Then the differential for intersection
homology vanishes.

\section{Appendix}
The most frequently applied corollary from the Decomposition Theorem
is \ref{cdt}. Its proof does not demand all the machinery of
\cite{BBD}. The Corollary \ref{cdt} can be deduced from functoriality
of intersection cohomology, which in turn relies only on a certain
vanishing statement called {\it Local Hard Lefschetz} in \cite{We2}.
We want to show how to construct inclusion and retraction in
\ref{cdt} using the induced maps constructed in \cite{BBFGK}.

We assume that $X$ and $Y$ are irreducible, $\dim(X)=n$ and
$\dim(Y)=m$. The map $f:X\rightarrow Y$ is proper and surjective.

\vskip6pt\noindent\hskip4pt{\sc Inclusion}: $\IC_Y\mpr{i}
Rf_*\IC_X$ is a map realizing $f^*:\IH^*(Y)\rightarrow\IH^*(X)$.
Such a map exists by \cite{BBFGK}

\vskip6pt\noindent\hskip4pt{\sc Retraction}:
$Rf_*\IC_X\mpr{r}\IC_Y$. Let $Z$ be a subvariety of $X$ such,
that $f_{|Z}$ is surjective and $\dim(Z)=\dim(Y)$. Such a construction has
already been applied in Remark \ref{ac} and \cite{BBFGK}
p.~173. We consider the sheaf $\IC_Z$ as a sheaf on $X$ supported by
$Z$. As in \S\ref{szesc} we construct a map $\alpha:\IC_Y\rightarrow
Rf_*\IC_Z$. We obtain a map of Verdier dual sheaves
$$Rf_*\IC_Z[2m]=Rf_*D\IC_Z=DRf_*\IC_Z
\mpr{D\alpha}D\IC_Y=\IC_Y[2m]\,.$$ Again, there exists a map
$\beta:\IC_X\rightarrow\IC_Z$ which agrees with the restriction of
the constant sheaf $\Cset_X\rightarrow\Cset_Z$. We compose
$Rf_*\beta$ with $D\alpha[-2n]$ and obtain a map
$r':Rf_*\IC_X\rightarrow\IC_Y$. Generically the composition
$r'\circ i$ is the multiplication by the degree of $f_{|Z}$. We
define $r={\rm deg}(f_{|Z})^{-1}r'$. Now, $r\circ
i:\IC_Y\rightarrow \IC_Y$ is generically identity. By \cite{Bo},
V.9.2 p.144 each morphism of intersection complexes is determined
by its restriction to an open set, therefore $r\circ
i=Id_{\IC_Y}$. \qed

\begin{rem1}\rm The existence of the induced maps constructed in
\cite{BBFGK} follows from local vanishing (Local Hard Lefschetz
\cite{We2}). Therefore these maps exist with rational
coefficients. Our inclusion and retraction are constructed for
rationals as well.\end{rem1}

\begin{rem1}\rm The assumption that we stay in the algebraic category
is essential. In the Example \ref{anal} there is a natural map
$\Cset_Y = \IC_Y\rightarrow Rf_*\IC_X = Rf_*\Cset_X$, but it is not
split injective. The construction of a retraction breaks down since
one cannot find a subvariety $Z\subset X$ with the desired
properties.\end{rem1}

Andrzej Weber

Instytut Matematyki, Uniwersytet Warszawski

ul. Banacha 2, 02--097 Warszawa, Poland

e-mail: aweber@mimuw.edu.pl

\end{document}